\documentclass[11pt]{article}
\usepackage{amsmath,amsfonts,amssymb,amsthm}
\newtheorem{theorem}{Theorem}[section]
\newtheorem{lemma}[theorem]{Lemma}
\newtheorem{corollary}[theorem]{Corollary}

\newtheorem{proposition}[theorem]{Proposition}
\newcommand{\R}{\mathbb{R}}            
\newcommand{\C}{\mathbb{C}}            

\newcommand{\Hb}{{\cal H}}
\newcommand{\del}{\partial}
\newcommand{\dba}{\overline{\partial}}

\newcommand{\Om}{\Omega}

\newcommand{\dinv}{\partial_{s}^{-1}}
\newcommand{\dnmap}{\Lambda_{b_1,b_2}}

\newcommand{\Skk}{{{\mathcal S}_k}}

\newcommand{\eps}{\varepsilon}
\newcommand{\dbar}{\bar{\partial}}

\newcommand{\pt}{{\tilde{p}}}
\newcommand{\rt}{{\tilde{r}}}

\newcommand{\ol}{\overline}
\newcommand{\dbinv}{{\dba^{-1}}}
\newcommand{\deinv}{{\del^{-1}}}
\DeclareMathOperator{\im}{Im}
\DeclareMathOperator{\re}{Re}

\title{On the Scattering Method for the
$\dba$-equation and Reconstruction of Convection Coefficients}

\begin{document}
\author{Alexandru Tamasan\thanks{This
work was done during the author's visit at IPAM-UCLA in the Fall
2003}\\
\parbox{12cm}{\small Department of Mathematics, University of
Toronto, 100 St. George St., Toronto, ON, M5S 3G4, Canada}}

\maketitle
\begin{center}
\parbox{12cm}{\small {\bf Abstract:}
In this paper we reconstruct convection coefficients from boundary
measurements. We reduce the Beals and Coifman formalism from a
linear first order system to a formalism for the $\dba$-equation.}
\end{center}

\section{Introduction}

The pioneering work of Nachman and Ablowitz
\cite{nachmanAblowitz}, Sylvester and Uhlmann
\cite{sylvesterUhlmann87}, Nachman \cite{nachman88}and Henkin and Novikov
\cite{henkinNovikov} introduced
inverse scattering methods to the parameter identification
problems. In their work, the linear Schr\"{o}dinger equation in
the physical space is paired with a pseudo-analytic equation in
the complex space of the parameter. Another method, due to Beals
and Coifman \cite{bealsCoifman}, pairs a first order $\dba$ system
in the physical space with a pseudo-analytic matrix equation in
the parameter space. Sung analyzed lower regularity assumptions in
\cite{sung1,sung2,sung3}. This method was ingeniously used by
Brown and Uhlmann \cite{brownUhlmann} in unique identification of
the conductivity $\sigma$ in $\nabla\cdot\sigma\nabla u=0$ and by
Cheng and Yamamoto \cite{chengYamamoto98}, \cite{chengYamamoto} in
proving unique
determination of the convection coefficients $b_1$ and $b_2$ in
$\Delta u + b_1u_x +b_2 u_y=0$.

We consider here the scattering problem for $\dba$- equations
(theorems 1.1 and 1.2 below). Here $\dba=(\del_x+i\del_y)/2$ is
the Cauchy-Riemann operator. This can be seen as a diagonal
version of the formalism in Beals and Coifman, see lemma
\ref{equivalence}. Due to the symmetry between the scattered
solutions in the physical space and the ones in the parameter
space, we are able to present a non-linear analog of the Fourier
inversion formula (compare \eqref{scatt_transf} and
\eqref{fourierInv} below).

As an application, we revisit the inverse problem proposed in 
\cite{chengYamamoto} and present a reconstruction
procedure. The method is based on
solving a singular boundary integral equation in the Hardy space
of functions in the exterior of the disc. This method was first
introduced by Knudsen and Tamasan in
connection with the electrical impedance tomography problem in
\cite{knudsenTamasan1}. The method presented here can be seen as its
generalization. 

I was informed recently about the reconstruction step being obtained 
independently by Tong, Cheng and Yamamoto \cite{tongCY}. 
I thank them for letting me know about their new result. The
main difference of the method presented here from their method is the 
formalism of in inverse scattering.

For $k\in\C$ arbitrarily fixed, we say that $u$ behaves like
$e^{izk}$ (written $u\sim e^{izk}$) in $L^r(\R^2_z)$ for large
$z$, if $u(z,k) e^{-izk}-1\in L^{r}(\R^2_z)$. We use the notation
$\langle k\rangle=(1+|k|^2)^{1/2}$. The scattering method is the content
of the following two theorems.
\begin{theorem}[Forward Scattering]\label{forward}
Assume that $q\in L^\pt_c(\R^2_z)$, $\pt >2$ has compact support.
For each $k\in \C$, the equation
\begin{equation}\label{phys_eq}
\frac{\del \Psi}{\del \ol{z}}(z)+q(z)\ol{\Psi}(z)=0,\,\,\, z\in
\C,
\end{equation}
has unique solutions $\Psi_r(z,k)\sim e^{izk}$ and
$\Psi_i(z,k)\sim ie^{izk}$ in $L^\pt$ for large $z$, and the
scattering transform
\begin{equation}\label{scatt_transf}
 t(k)=-\frac{i}{\pi}\int_{\R^2}e^{i\ol{zk}}\ol{q}(z)\left(
\Psi_r(z,k)-i\Psi_i(z,k)\right)d\mu(z)
\end{equation}is well defined.
Moreover, if $q\in W^{\eps,\pt}_c(\R^2_z)$ for some $\eps>0$ and
$k \in \C-\{0\}$, we have
\begin{align}\label{k_estimates}
\|\Psi_r(z,k)e^{-izk}-1\|_{L^\pt
(\R^2_z)}+\|\Psi_i(z,k)]e^{-izk}-i \|_{L^\pt (\R^2_z)}\leq C
\langle k\rangle^{-\eps}
\end{align}and
\begin{align}
\|[\Psi_r(z,k)-i\Psi_i(z,k)]e^{-izk}-2 \|_{W^{1,\pt} (\R^2_z)}\leq
C \langle k\rangle^{-\eps},\label{extrasmooth}
\end{align}
and then $t\in L^r(\R^2_k)$ for each $r>2/(\eps +1)$.  In
particular $t\in L^r(\R^2_k)\cap L^{r'}(\R^2_k)\cap
L^{\rt}(\R^2_k)$ for some $r<2$, where $\rt^{-1}=r^{-1}-1/2$ and $
r'^{-1}+r^{-1}=1$.
\end{theorem}

\begin{theorem}[Inverse Scattering]\label{inverse} Let $q$, $\Psi_r$, $\Psi_i$ and $t(k)$ and $r$, $r'$, $\rt$
be as given in the forward scattering. Then the equation
\begin{equation}\label{freq_eq}
\frac{\del \Phi}{\del \ol{k}}(k)+t(k)\ol{\Phi}(k)=0,\,\, k\in \C,
\end{equation}
has unique solutions $\Phi_r\sim e^{izk}$ and $\Phi_i\sim
ie^{izk}$ in $L^\rt(\R^2_k)$ for large $k\in\C$. Moreover,
$\Psi$'s and $\Phi$'s are related by
\begin{align}\label{identities}
\re{\Phi_i}=-\im{\Psi_r},\qquad& \re{\Phi_r}=\re{\Psi_r},\\
\im{\Phi_i}=\im{\Psi_i}, \qquad
&\im{\Phi_r}=-\re{\Psi_i},\nonumber
\end{align}in particular $
\Phi_r-i\Phi_i=\Psi_r-i\Psi_i$ and
\begin{align}\label{fourierInv}
q(z)=-\frac{i}{\pi}\int_{\R^2}e^{i\ol{zk}}\ol{t}(k)\left(
\Phi_r(z,k)-i\Phi_i(z,k)\right)d\mu(k).
\end{align}
\end{theorem}

Let $\Om\subset\R^2$ be a bounded, simple
connected domain with Lipschitz boundary and $\pt>2$. For
$b_1,b_2\in L^\pt(\Om)$ and $g\in W^{2-1/\pt,\pt}(\del\Om)$, let
$u\in W^{2,\pt}$ be the unique solution of the boundary value
problem
\begin{align}\label{conveq}
&\Delta u(x)+b_1(x)\frac{\del u}{\del x_1}(x) + b_2(x)\frac{\del u}{\del x_2}(x) =0, \qquad x\in \Om\\
&u(x)= g(x),\qquad x\in\del\Om\nonumber.
\end{align} The Dirichlet to Neumann map $\dnmap:W^{2-1/\pt,\pt}(\del\Om)\to
W^{1-1/\pt,\pt}(\del\Om)$ is given by $$\dnmap g(x)=
\nu_1(x)\frac{\del u}{\del x_1}(x) + \nu_2(x)\frac{\del u}{\del
x_2}(x),\qquad x\in\del\Om,$$where $(\nu_1(x),\nu_2(x))$ is the
outer normal at $x$ on the boundary. Cheng and Yamamoto proved
that $\dnmap$ uniquely determines $b_1$ and $b_2$ in $L^\pt(\Om)$.

Working with the equation in the whole plane and using the inverse
scattering for $\dba$-equations allows us to go beyond uniqueness
and present a method of reconstruction.  We prove
the following result.
\begin{theorem} Let $\Om\subset\R^2$ be
bounded, simple connected domain with Lipschitz boundary, and
$b_1,b_2\in W^{\eps,\pt}_c(\Om)$, with support inside $\Om$ for
some $\eps >0$ . Then $b_1,b_2$ can be reconstructed from
$\dnmap$.
\end{theorem}The fact that they vanish on the boundary is not a severe
restriction as one can always extended the coefficients across the
boundary, preserving the regularity, and then have them vanish
outside a ball. The Dirichlet-to-Neumann map can be pushed to an
outside boundary as shown by Nachman in \cite{nachman88}, see also
\cite{knudsenTamasan}. While $L^\pt(\Om)$ is enough regularity to
prove unique determination of $b_1,b_2$, we assume here
$\epsilon$-extra regularity and provide a reconstruction method.

In the end we point out the connection with the first order $\dba$
system and characterize its Cauchy data in terms of the
Dirichlet-to-Neumann map of a related second order elliptic
equation, thus answering a question of Uhlmann in \cite{uhlmann}.

\section{Proof of the theorems \ref{forward} and \ref{inverse} }
We identify a point in $\R^2$ with a point in the complex plane by
$x_1+ix_2=z$. By $\dbinv$ we denote the solid Cauchy transform
\begin{align}
\dbinv
f(z)=\frac{1}{\pi}\int_{\R^2}\frac{f(\zeta)}{z-\zeta}d\mu(\zeta),
\end{align}where $d\mu(\zeta)$ is the Lebesgue area.
We also denote by $e(z,k)=\exp(i(zk+\ol{zk}))$.

We look for solutions of \eqref{phys_eq} of the form
$\Psi_r=\psi_r e^{izk}$ and $\Psi_i=i\psi_i e^{izk}$ with $\psi_r,
\psi_i\in 1+L^\pt(\R^2_z)$. The equations for $\psi_r$
respectively $\psi_i$ are
\begin{align}\label{eqpsis}
\frac{\del}{\del \ol{z}}\psi_r+qe(z,-k)\ol{\psi_r}&=0,\\
\frac{\del}{\del \ol{z}}\psi_i-qe(z,-k)\ol{\psi_i}&=0.\nonumber
\end{align}
A key ingredient is the Hardy-Littlewood Sobolev inequality which
yields $\dba^{-1}:L^p(\R^2)\to L^\pt(\R^2)$ is bounded (see Stein
\cite{stein}) for $p$ and $\pt$ related by
\begin{align}\label{tilde}
\frac{1}{\pt}=\frac{1}{p}-\frac{1}{2}.
\end{align}
Since $q\in L^\pt_c(\R^2)\subset L^2(\R^2)$ and $L^\pt(\R^2) \cdot
L^2(\R^2)\subset L^p(\R^2)$ we have $\dbinv(q\cdot): L^pt(\R^2)\to
L^\pt(\R^2)$ is a bounded operator. Since $q$ has compact support
we can use Rellich imbedding to conclude that $\dbinv(q\cdot):
L^p(\R^2)\to L^\pt(\R^2)$ is compact. Then we can apply Fredholm's
alternative in $L^\pt(\R^2)$ to the equivalent integral equation
\begin{align}\label{integral}
\{I +\dba^{-1} [q(\cdot)e(\cdot,-k)\ol{(\cdot)}]\}(\psi_r(z)-1) &=
\dba^{-1}[qe(\cdot,-k)\ol{}],\\
\{I -\dba^{-1} [q(\cdot)e(\cdot,-k)\ol{(\cdot)}]\}(\psi_i(z)-1)
&=- \dba^{-1}[qe(\cdot,-k)\ol{}].\nonumber
\end{align}The fact that the homogeneous equation has only the null solution
comes from Liouville's theorem for pseudo-analytic functions with
coefficients in $L^\pt(\R^2_z)\cap L^p(\R^2_z)$ shown by Vekua
\cite{vekua}. Since we integrate in \eqref{scatt_transf} over the
support of $q$, together with the imbedding $L^\pt_{loc} \subset
L^p_{loc}$, gives a pointwise well defined $t(k)$.

For $k\in\C-0$ let $(\dba-i\ol{k})^{-1}$ be defined by
$e(z,k)\dbinv(e(z,-k)\cdot)$ and let the indexes $\pt$ and $p$ be
related by \eqref{tilde}. An interpolation (with $\epsilon$ being
the interpolation parameter) between the estimates of Nachman
\cite{nachman96} $||(\dba-i\ol{k})^{-1} f||_{L^\pt}\leq
C||f||_{L^p(\R^2)}$ and  $||(\dba-i\ol{k})^{-1} f||_{L^\pt}\leq
(C/|k|)||f||_{W^{1,p}(\R^2)}$ gives
\begin{align}\label{estimate}
||(\dba-i\ol{k})^{-1}
f||_{L^\pt}\leq\frac{C}{|k|^{\epsilon}}||f||_{W^{\epsilon,p}(\R^2)}.
\end{align}
See Proposition 2.3 in \cite{knudsenTamasan} for details. Since
$|e(z,k)|=1$ the last estimate implies that
\begin{align*}
||\dbinv(e(\cdot, -k)q)||_{L^\pt(\R^2_z)}\leq C\langle
k\rangle^{-\epsilon}||q||_{W^{\eps,p}(\R^2_z)}.
\end{align*}The decay rate in \eqref{k_estimates}
follows from the (uniform in $k$) bounded-ness of the map
$[I-\dba^{-1}qe(\cdot,-k)\cdot]^{-1}$ from $L^\pt(\R^2_z)$ to
$L^\pt(\R^2_z)$ as explained above. The further regularity
property for the combination $\Psi_r+i\Psi_i$ in
\eqref{extrasmooth} will be shown in Lemma \ref{equivalence}. For
now we assumed it holds.

Brown and Uhlmann \cite{brownUhlmann} showed that $q\in
L^\pt_c(\R^2_z)$ implies $t\in L^2(\R^2_k)$. While this is good
enough for existence, for reconstruction we need $t\in L^r$ for
some index $r<2$. This is ensured by extra regularity imposed in
$q$ as was shown by Knudsen and the author in
\cite{knudsenTamasan}. For completeness we repeat the arguments.
The main ingredient is an $L^2$ bounded-ness property for
pseudo-differential operators with non smooth symbol (see Coifman
and Meyer \cite{meyer} or Brown and Uhlmann \cite{brownUhlmann2}).
If $Mq$ is defined by
$$
Mq(k)=\frac{i}{\pi}\int_{\R^2}e(z,k)\ol{q}(z)a(z,k)d\mu(z),
$$ where $a$ has compact support in $z$ and $||a(\cdot,k)||_{H^1(\R^2_z)}\leq
{C}\langle k\rangle^{-\epsilon}$, then $M: L^2(\R^2_z)\to
L^2(\R^2_k)$ is bounded.

Rewrite now
\begin{align}
 t(k)&=-2i {\cal F}(\ol{q})(-2k_1,2k_2)+T(k)\\
 T(k)&=\frac{i}{\pi}\int_{\R^2}e(z,k)\ol{q}(z)\left[
\psi_r(z,k)+\psi_i(z,k)-2\right]d\mu(z),
\end{align}
where ${\cal F}$ is the Fourier transform. Since $q\in
L_c^\pt(\R^2_z)\subset L^r(\R^2_z)$ for $1\leq r\leq 2$ we get
${\cal F}(\ol{q})\in L^s(\R^2_k)$ for $s>2/(1+\epsilon)$.

Let $M$ be the operator defined by
$a(z,k)=\chi(z)[\psi_r(z,k)+\psi_i(z,k)-2]\in
W^{1,\pt}_c(\R^2_z)\subset H^1(\R^2)$, where $\chi$ is a cut-off
function equal to $1$ on the support of $q$. The following chain
of inequalities for $0<\delta<\epsilon$ give the result
$$\|T\|_{L^s(\R^2_k)}=\|\langle k\rangle^{-\delta}M\ol{q}\|_{L^s(\R^2_k)}
\leq C||\langle
k\rangle^{-\delta}||_{L^{(1/s-1/2)^{-1}}}\|q\|_{L^2(\R^2_z)},$$
for $\delta >2(1/s-1/2)$ or equivalently
$s>2/(\delta+1)>2/(\epsilon+1)$.

In order to exhibit the relation with the old formalism, we prove
theorem \ref{inverse} by reducing it to the former. Let us define
$m_1(z,k)$ and $m_2(z,k)$ in terms of the $\psi$'s by
\begin{align}\label{tranzition}
m_1(z,k)&=\frac{1}{2}\left(\psi_r(z,k)+\psi_i(z,k)\right)\\
m_2(z,k)&=\frac{1}{2}e(z,-k)\left(\ol{\psi_i}(z,k)-\ol{\psi_r}(z,k)\right)\nonumber
\end{align}
The simple result below shows that $(m_1,m_2)^t$ is the first
column of the Jost matrix in the complex geometrical optic
solutions of Beals and Coifman.
\begin{lemma}\label{equivalence}Let $m_1$ and $m_2$ defined in
\eqref{tranzition}. Then $m_1(\cdot,k) -1,m_2(\cdot,k)\in
L^\pt(\R_z)$, and they satisfy
\begin{align}\label{Dsys}
\dba m_1&= q m_2\\
(\del+ik)m_2&= \ol{q}m_1.\nonumber
\end{align}
Moreover, the following estimates hold,
\begin{align}
||m_1(\cdot,k)-1||_{W^{1,\pt}(\R^2_z)}\leq C\langle k\rangle^{-\epsilon}\\
||m_2(\cdot,k)||_{L^\pt(\R^2_z)}\leq C\langle
k\rangle^{-\epsilon}.
\end{align}
\end{lemma}
\begin{proof}From their definition $m_1(\cdot,k)-1, m_2(\cdot,k)\in
L^\pt(\R^2_z)$ since $\psi_r,\psi_i\in L^\pt(\R^2_z)$ and
$|e(z,k)|=1$. The fact that they solve the system \eqref{Dsys}
comes from a straightforward calculation and the equations
\eqref{eqpsis}. The $L^{\pt}$ estimates of decay in $k$ for noth
$m_1$ and $m_a$ come from the estimates \eqref{k_estimates} for
$\psi_r$ and $\psi_i$ proven above. We are left to justify the
extra smoothness gained by $m_1$. From the first equation we have
that $m_1-1=\dbinv (qm_2)$. Since $q\in
 L^\pt_c(\R^2)\subset L^2(\R^2)$ and $L^2\cdot L^\pt\subset L^p$ we
have $\dbinv (qm_2)\in W^{1,\pt}(\R^2_z)$ with an imbedding
constant which depends on the support of $q$ but it is independent
of $k$. We have the following chain of inequalities.
\begin{align*}
||m_1(\cdot,k)-1||_{W^{1,\pt}(\R^2_z)} &= \|\dbinv
(qm_2)\|_{W^{1,\pt}(\R^2_z)}\leq C\| q(\cdot)
m_2(\cdot,k)\|_{L^p(\R^2_z)}\\
&\leq C\| q\|_{L^2(\R^2)}\| m_2(\cdot,k)\|_{L^\pt(\R^2_z)}\leq
C<k>^{-\epsilon}.
\end{align*}
\end{proof}
This also completes the proof the theorem \ref{forward}.

Formulate the inverse scattering formalism of Beals and Coifman
only in terms of the first column of Jost matrix, see Knudsen and
Tamasan \cite{knudsenTamasan} for details. For the analysis with
$q\in L^1(\R^2)\cap L^\infty(\R^2)$ see Sung \cite{sung1}, or
Brown and Uhlmann \cite{brownUhlmann} for $q\in L^\pt_c(\R^2)$.

\begin{theorem}[Beals \& Coifman scattering method]
Let $q \in W^{\epsilon ,\pt}_c(\R^2)$, For any $z\in C$, the
system \eqref{Dsys} has a unique solution $m_1(z,k)$, $m_2(z,k)$
with $ (m_1(\cdot,k) - 1, m_2(\cdot,k))\in L^\pt(\R^2_z)$.
Furthermore, the map $k\to m(\cdot ,k)$ is differentiable (in the
norm topology) with values in $W^{\epsilon,\pt}(dx)$ and satisfies
pointwise in $z\in \C$ the system
\begin{align}\label{dbar}
\frac{\del}{\del \ol k}m_1(z,k) = t(k) e(z,-k) \ol{m_2(z,k)},\\
\frac{\del}{\del \ol k}m_2(z,k) = t(k) e(z,-k)
\ol{m_1(z,k)},\nonumber
\end{align}
where
\begin{equation}\label{S}
  t(k)= -\frac{i}{\pi}\int_{\R^2}e(z,k)\ol q(z) m_{1}(z,k)d\mu(z).
\end{equation}
\end{theorem}

Look for solutions of \eqref{freq_eq} in the form $\Phi(z,k)=i
e^{izk}\phi_r(z,k)$ respectively $\Phi_i=e^{izk}\phi_r(z,k)$. As
in the forward problem, they must satisfy an integral formulation
analogous to \eqref{integral} where the r\^{o}le of $k$ and $z$ is
reversed. Since $t(k)\in L^r(\R^2_k)\cap L^2(\R^2_k)$ we have
existence and uniqueness for their solution in $L^\rt(\R^2_k)$,
where $\rt^{-1}=r^{-1}-1/2$. Using the equations \eqref{dbar} it
is easy to check that
\begin{align}
\frac{\del}{\del\ol{k}}(m_1-m_2)(z,k)&=-t(k)e(z,-k)\ol{m_1-m_2}(z,k)\\
\frac{\del}{\del\ol{k}}(m_1+m_2)(z,k)&=
t(k)e(z,-k)\ol{m_1+m_2}(z,k).\nonumber
\end{align}
By the uniqueness result we must have
\begin{align}
\phi_i(z,k)&=m_1(z,k)+m_2(z,k)\\
\phi_r(z,k)&=m_1(z,k)-m_2(z,k).\nonumber
\end{align}The following equalities show the relation between
solutions of the forward and inverse equation.
\begin{align*}
\Phi_i&=ie^{izk}\phi_i=
ie^{izk}(m_1+m_2)=\frac{ie^{izk}}{2}(\psi_r+\psi_i)+
\frac{ie^{izk}}{2}e(z,-k)(\ol{\psi_i}-\ol{\psi_r})\\
&=\frac{i}{2}\Psi_r +\frac{1}{2}\Psi_i
-\frac{1}{2}\ol{\Psi_i}-\frac{i}{2}\ol{\Psi_r}=-\im{\Psi_r}+i\im{\Psi_i}.
\end{align*}
Similarly, $\Phi_r=\re{\Psi_r}-i\re{\Psi_i}$. These prove the
identities \eqref{identities}. Formula \eqref{fourierInv} is due
to a symmetry argument as follows. Starting with $q$ produce
$\psi_r$ and $\psi_i$ by solving \eqref{phys_eq}. Via
\eqref{identities} produce $\phi_r$ and $\phi_i$ and then $t(k)$
as in \eqref{fourierInv}. Take this $t(k)$ and do now forward
scattering starting from the $k$-space, i.e. produce $\Phi_r$ and
$\Phi_s$ by solving \eqref{freq_eq} and via \eqref{identities}
produce $\Psi_r$ and $\Psi_i$. Define a potential $q_1$ using
\eqref{fourierInv} for the $z$-space. In particular we know that
for any $k\in\C$ we have $\dba\Psi_r +q\ol{\Psi_r}=0$ since we
started that way, but also now we have $\dba\Psi_r
+q_1\ol{\Psi_r}=0$. In particular we have
$(q(z)-q_1(z))\ol{\Psi_r}(z,k)=0$ for all $k\in\C$. Hence $q=q_1$.

\section{Reconstructing convection coefficients}
In this section we apply the above scattering method to
reconstruction of the convection coefficients $b_1,b_2$ in
\begin{align}\label{elliptic}
\Delta u(x)+b_1\frac{\del u}{\del x}(x)+ b_2(x)\frac{\del u}{\del
x}(x)=0, x\in \Om
\end{align}
from the Dirichlet-to-Neumann map $\dnmap$. Here $\Om\subset\R^2$
is a bounded, simply connected domain with Lipschitz boundary.

We assume here that $b_1,b_2\in W^{\epsilon,\pt}_c(\Om)$, $\pt>2$
are real valued maps with compact support in $\Om$ and set
$b=(b_1+ib_2)/4$.

The following result from Vekua \cite{vekua} makes the reduction
of \eqref{elliptic} to a $\dba$-equation. If $u$ is a solution of
\eqref{elliptic} then $w=\del u$ solves
\begin{align}\label{dbarlong}
\dba w(z) + \ol{b}(z)\ol{w}(z) +b(z)w(z) =0.
\end{align}
\begin{lemma}\label{reduction} Let $\Om$ be simply connected with Lipschitz boundary. If $u\in W^{2,\pt}(\Om)$ is a solution of
\eqref{elliptic}, then $w=\del u\in W^{1,\pt}(\Om)$ is a solution
of \eqref{dbarlong}. Conversely, if $w\in W^{1,\pt}(\Om)$ is a
solution of \eqref{dbarlong} then there exists an $u\in
W^{2,\pt}(\Om)$ solution of \eqref{elliptic} and such that $\del
u=w$ in $\Om$.
\end{lemma}

\begin{proof} By Sobolev imbedding we have $u\in C^{1+\alpha}(\Om)$ with $\alpha
=1-1/\pt$ and $w\in C^{\alpha}(\Om)$. As a direct consequence of
the Poincar\'e lemma, notice that if $\dbar w$ is real valued,
then $w= \del u$ for some real valued $u$. Indeed $2\dbar w
=(\del_x +i\del_y)(f+ig)=(\del_x f-\del_y g) +i(\del_x g+\del_y
f).$ By assumption $\del_x g = -\del_y f$, from where the one-form
$gdy-fdx$ is exact. Therefore, there exists a real valued $F$ such
that $dF=(-f)dx +gdy$. We have $w= f +ig =\del_x (-u)
-i\del_y(-u)= \del (-2F)$. The equivalence is now apparent.
\end{proof}

Now we extend $b\in W^{\eps,\pt}_c(\Om)$ by zero outside $\Om$.
Its extension denoted also by $b$ preserves regularity $b\in
W^{\eps,\pt}_c(\R^2)$. From now on we shall work with solutions of
\eqref{dbarlong} in the whole plane.
\begin{lemma}\label{Ws} The equation \eqref{dbarlong} has unique solutions
in the whole plane  $W_r(z,k)\sim e^{izk}$ respectively
$W_i(z,k)\sim ie^{izk}$ in $L^\pt(\R^2_z)$ for large $z$.
Moreover, $e^{-izk}W_r-1,e^{-izk}W_i-i\in W^{1,\pt}(\R^2_z)$ and
$W_r(\cdot,k),W_i(\cdot,k)\in W^{1,\pt}_{loc}(\R^2)$.
\end{lemma}
\begin{proof}
As in the proof of theorem 1.1, we look for solutions
$W(z,k)=e^{izk}w(z,k)$ with $w-1 \in L^\pt$. The equation for $w$
is
\begin{align}
\dba (w(z)-1) + b(z)(w(z)-1) + e(z,-k)\ol{b}(z)(\ol{w}(z)-1) =
-b(z)-e(z,-k)\ol{b}.
\end{align}Using the fact that $\dbinv :f\in L^\pt_c(\R^2)\mapsto W^{1,\pt}(\R^2)$
together with $b$ of compact support we get $\dba^{-1}
(b\cdot):L^\pt(\R^2_z)\to L^\pt(\R^2_z)$ is a compact operator. We
apply Fredholm's alternative in $L^\pt(\R^2_z)$ to the equivalent
integral equation
$$\{[I +\dba^{-1} [b(\cdot)+e(\cdot,-k)\ol{b}\ol{(\cdot)}]\}(w(z)-1) =
-\dba^{-1}[b+e(\cdot,-k)\ol{b}].$$ Uniqueness comes from
Liouville's theorem for the $\dba$-equation with coefficients in
$L^\pt(\R^2)\cap L^p(\R^2)$, see Vekua \cite{vekua}. By
construction we already have that $g=w_r-1\in W^{1,\pt}(\R^2_z)$.
Then $W_r(z,k)=e^{izk}(g+1)\in L^\pt_{loc}(\R^2_z)$, $\del
W_r=ike^{izk}g(z,k)+ike^{izk} + e^{izk}\del g\in
L^\pt_{loc}(\R^2_z)$ and $\dba W_r= e^{izk}\dba g\in
L^\pt_{loc}(\R^2_z)$. Similar relations hold for $W_i$.
\end{proof}

To simplify notations, let
\begin{align}\label{new_q}
q(z)=\ol{b}(z)e^{\dbinv b (z) -\deinv\ol{b}(z)}
\end{align}denote a new potential and notice that if $w$ is a solution of
\eqref{dbarlong} then $v=e^{\dbinv b} w$ is
a solution of
\begin{align}\label{dbareq}
\dba v +q\ol{v}=0.
\end{align}
Since $b\in L^\pt(\R^2)\cap L^p(\R^2)$ we have that $\dbinv b\in
L^\infty(\R^2)\cap C^\frac{r-2}{2}$, see Vekua \cite{vekua}. Then
$e^{-\dbinv b}\in L^\infty(\R^2)$ and so $q\in L^\pt(\R^2)\cap
L^p(\R^2)$ .

The next theorem relates scattering solutions of \eqref{dbarlong}
to scattering solutions of \eqref{dbareq} and gives the behavior
in $k$ of $W_r(z,k)$ and $W_i(z,k)$.
\begin{proposition}\label{relation} Let $b\in W_c^{\eps,\pt}(\R^2)$, for some $\eps >0$. Let $W_r$
and $W_i$ be the scattering solutions for \eqref{dbarlong} as
given by the lemma above, and let $\Psi_r$ and $\Psi_i$ be the
scattering solutions of \eqref{dbareq} as given by the theorem
1.1. Then $W_r= {e^{-\dbinv b}} \Psi_r$, $W_i=e^{-\dbinv b}
\Psi_i$ and
\begin{align}\label{longKestimate}
\|W_r(z,k)e^{-izk}- e^{-\dbinv b}\|_{L^\pt(\R^2_z)}+
\|W_i(z,k)e^{-izk}-ie^{-\dbinv b}\|_{L^\pt(\R^2_z)}\leq C\langle k\rangle^{-\eps},\\
\left\|[W_r(z,k)-iW_i(z,k)]e^{-izk}- 2e^{-\dbinv
b}\right\|_{W^{1,\pt}(\R^2_z)}\leq C\langle
k\rangle^{-\eps}.\nonumber
\end{align}
\end{proposition}
\begin{proof}The fact that $W_r$ and $W_i$ solve \eqref{dbarlong} is
trivial. Uniqueness result of lemma \ref{Ws} ensures that they are
the scattering solutions of \eqref{dbarlong}. The estimates follow
directly from the estimates for $\psi_r$ and $\psi_i$ in
\eqref{k_estimates} and \eqref{extrasmooth} and from the fact that
$e^{-\dbinv b}\in L^\infty(\R^2_z)$ as noticed before. Again, the
imbedding $W^{1,\pt}(\R^2) \subset C^{1-2/\pt}(\R^2)$ shows that
the estimates \eqref{longKestimate} hold pointwise in $z\in\C$.
\end{proof}

We have now all the ingredients necessary for reconstruction.
Since $q$ in \eqref{new_q} has compact support in $\Om$, the
scattering transform depends only on the traces on $\del\Om$ of
the scattering solutions $\Psi_r$ and $\Psi_i$. Let
$\nu=\nu_1+i\nu_2$ be the complex-normal to the boundary. Then
\begin{align}\label{newScatT}
t(k)&=-\frac{i}{\pi}\int_{\Om}e^{i\ol{zk}}\ol{q}(z)\left(
\Psi_r(z,k)-i\Psi_i(z,k)\right)d\mu(z)
=\frac{i}{\pi}\int_{\Om}e^{i\ol{zk}}\left(
\del\ol{\Psi_r}(z,k)-i\del\ol{\Psi_i}(z,k)\right)d\mu(z)\nonumber\\
&=\frac{i}{2\pi}\int_{\del\Om}e^{i\ol{zk}}\ol{\nu}(z)
\left(\ol{\Psi_r}(z,k)-i\ol{\Psi_i}(z,k)\right)d\sigma(z),
\end{align} The last equality uses the fact that
$\del(e^{i\ol{zk}})=0$ .

Next we show how to reconstruct traces of $\Psi_r$ and $\Psi_i$ to
$\del\Om$ from the Dirichlet to Neumann map $\dnmap$. First we
reconstruct traces of $W_r$ and $W_i$ to $\del\Om$.

As in Knudsen and Tamasan \cite{knudsenTamasan}, we consider the
single layer potential operator $\Skk :C^\alpha(\del\Om)\to
C^\alpha(\del\Om)$, $\alpha=1-2/\pt$, defined by
$$\ol{\Skk} f(z)=\frac{1}{2\pi i}p.v.\int_{\del\Om}{f(\zeta)}{\ol{g_k}(\zeta-z)}d\ol{\zeta},\qquad z\in\del\Om,$$
where $g_k(z)=e^{-izk}/(\pi z)$ is a Cauchy kernel for $\dba$
which also takes into account the exponential behavior at
infinity. For Lipschitz boundary $\Skk$ is a bounded operator
(e.g. see Muskhelishvili \cite{musk}). Since $q$ has compact
support we have that $W_r$ and $W_i$ are analytic outside $\Om$
and behaves like $e^{izk}$ at infinity. Traces of such functions
will satisfy a singular boundary equations involving $\Skk$.
Inside $\Om$ they satisfy a pseudo-analytic equation. This will
impose constrains (in terms of $\dnmap$) on their trace. We will
prove that these two conditions are sufficient to determine the
traces.

We notices already that $W_r(\cdot,k), W_i(\cdot,k)\in
C^\alpha(\R^2)$ with $\alpha=1-2/\pt$, whence their traces on
$\del\Om$ are in $C^\alpha(\del\Om)$. Let
$$C^\alpha_0(\del\Om):=\{h\in C^\alpha(\del\Om):\,\,\int_{\del\Om}h(s)ds=0
\}.$$ Define now a right inverse of the tangential vector field
$\del_s$ (here $s$ is the arc length) on $\del\Om$ by
\begin{align}
\dinv f(t)=\int_{0}^t f(s)ds,
\end{align}
for $f\in C^\alpha(\del\Om)$. In the above integral we fixed an
arbitrary point on $\del\Om$ from where we measure the arc length
counter-clockwise. Notice that $\dinv: C^\alpha_0(\del\Om)\to
C^{1+\alpha}(\del\Om)$ is a well defined (independent of the
reference point) bounded operator. The following result defines a
Hilbert transform for the pseudo-analytic maps.
\begin{lemma}$\Hb_b\equiv-\dnmap\dinv :C^\alpha_0(\del\Om)\to
C^{\alpha}(\del\Om)$ is a bounded operator.
\end{lemma}
\begin{proof} Let $g=\dinv f\in C^{\alpha+1}(\del\Om)\subset W^{2-1/\pt,\pt}$.
Classical theory of PDE (e.g. see Gilbarg and Trudinger
\cite{gilbarg}) gives that the boundary value problem
\begin{align}
&\Delta u(x)+b_1\frac{\del u}{\del x}(x)+ b_2(x)\frac{\del u}{\del
x}(x)=0,\,\,\, x\in \Om\\
&u|_{\del\Om}(x)=g(x),\,\,\, x\in\del\Om \nonumber
\end{align} has a unique solution up to a constant in
$W^{2,\pt}(\Om)$ and $||u||_{W^{2,\pt}(\Om)}\leq
C||g||_{W^{2-1/\pt,\pt}(\del\Om)}$. Using the mapping properties
of the Dirichlet to Neumann map we have
\begin{align*}
\|\Hb_b f\|_{C^\alpha(\del\Om)}&\leq\|\dnmap
g\|_{W^{1,\pt}(\del\Om)}\leq\|\nabla u\|_{W^{1,\pt}(\Om)}\leq
||u||_{W^{2,\pt}(\Om)}\\
&\leq C||g||_{W^{2-1/\pt,\pt}(\del\Om)}\leq
C||g||_{C^{1+\alpha}(\del\Om)}\leq C\|f\|_{C^\alpha(\del\Om)}.
\end{align*}
\end{proof}
Next we show that $\Hb_b$ reconstructs traces of the exponentially
growing solutions on $\del\Om$.
\begin{theorem}[Trace theorem]Let $b\in W_c^{\eps,\pt}(\Om)$.
Consider the class of functions $${\cal B}=\{h\in
C^\alpha(\del\Om):\im(\nu h)\in C^\alpha_0(\del\Om)\}.$$ Then, for
each $k\in \C$ arbitrarily fixed, the traces
$h_r=W_r(\cdot,k)|_{\del\Om}.$, respectively
$h_i=W_i(\cdot,k)|_{\del\Om}$ are the unique solution in ${\cal
B}$ of the systems
\begin{align}
&(I-i{\Skk})h_r(z)=2e^{izk},\qquad z\in\del\Om,\label{outside}\\
&{\cal H}_b(\im(\nu h_r))(z)=\re(\nu h_r)(z),\qquad
z\in\del\Om,\label{inside}
\end{align}respectively,
\begin{align}\nonumber
&(I-i{\Skk})h_i(z)= 2ie^{izk},\qquad z\in\del\Om,\\
&{\cal H}_b(\im(\nu h_i))(z)=\re(\nu h_i)(z),\qquad
z\in\del\Om.\nonumber
\end{align}
\end{theorem}
\begin{proof}
We argue only for $W_r$, the arguments for $W_i$ are similar.

We prove first the necessity. The arguments for \eqref{outside}
are identical to the ones in \cite{knudsenTamasan} reason for
which we only sketch them here. Fix $k\in \C$ and suppress the $k$
dependence, we have $W_r(\cdot) = W_r(\cdot,k)$ is analytic
outside $\Om$ and $e^{-izk} W_r-1\in L^\pt (\C-\Om)$. The
Green-Gauss formula for $z\in\C-\Om$ gives
\begin{align}\label{outcont}
e^{-izk} W_r(z)-1=-\frac{1}{2\pi i}\int_{\del\Om}\frac{e^{-i\zeta
k}W_r(\zeta)}{\zeta-z}d\zeta.
\end{align}
Now let $z$ approach (from the exterior) a boundary point $z_0$
and use Plemelj formula (see Muskhelishvili \cite{musk}).
$$\lim_{z\to z_0}\frac{1}{2\pi
i}\int_{\del\Om}\frac{f(\zeta)d\zeta}{\zeta-z}=-\frac{1}{2}f(z_0)+
\frac{1}{2\pi i}p.v.\int_{\del\Om}\frac{d\zeta}{\zeta-z_0}$$ to
get \eqref{outside}.

Next we prove the necessity of \eqref{inside}. Recall from lemma
\ref{Ws} that $W_r(z)=\del u(z)$ for some $u\in W^{2,\pt}(\Om)$
which solve the equation \eqref{elliptic}. Therefore
\begin{align}
h_r=W_r|_{\del\Om}=\frac{1}{2}(\del_x +i\del_y)u|_{\del\Om}.
\end{align}
For $z\in\del\Om$ let $(\nu_1(z),\nu_2(z))$ be the unit outer
normal, we also let $\nu(z)=\nu_1(z)+i\nu_2(z)$. Next we express
the partial derivatives for points on the boundary in terms of the
tangent $\del_s$ and the normal $\del_\nu$ derivatives
\begin{align}
\nabla u(x)=
\begin{pmatrix}
-\nu_2&\nu_1\\
\nu_1&\nu_2
\end{pmatrix}
\begin{pmatrix}
\del_s u\\
\dnmap u
\end{pmatrix},
\end{align}
where we recall $\del_\nu u=\dnmap u$.
Therefore $2h_r=(\del_x -i\del_y)u=-i\ol{\nu}\del_s u
+\ol{\nu}\dnmap u$, or, using $\nu\ol{\nu}=1$,
\begin{align}\label{brelation}
2\nu h_r=\dnmap u -i\del_s u.
\end{align}
Note that $\im (\nu h_r)=-\del_su/2$ and thus $h_r\in {\cal B}$
and $\dinv (\im(\nu h_r))$ makes perfect sense. Identifying the
real part in \eqref{brelation} gives \eqref{inside}. Notice not
only that we proved necessity but also we provided existence of
solutions for \eqref{inside} and \eqref{outside}.

Conversely, let $h\in{\cal B}$ be a solution of the system
\eqref{inside} and \eqref{outside}. We extend $h$ inside $\Om$ by
the following procedure. Inspired by \eqref{brelation} define
$g=-\dinv \im(2\nu h)\in C^{\alpha}(\del\Om)$ then uniquely solve
the boundary value problem \eqref{conveq} for $u\in
W^{2,\pt}(\Om)$. Notice $g$ is real valued hence $u$ has also real
values. Define $W_r(z) =\del u(z)$ inside $\Om$ and notice that
$\del u|_{\del\Om}\in C^\alpha(\del\Om)$. Now check that $\del
u|_{\del\Om} =h$. Indeed, as before,$ 2\del u=-i\ol{\nu}\del_s u
+\ol{\nu}\dnmap u= i\ol{\nu}\im(2\nu h)+\ol{\nu}\re(2\nu h)$. The
last equality used the fact that $h$ is a solution of
\eqref{inside}. Multiplication by $\nu$ gives $\del u= h$.

Inspired by \eqref{outcont} define $W_r$ analytically outside
$\Om$ by
\begin{align}
W_r(z)=e^{izk}-\frac{1}{2\pi i}\int_{\del\Om}\frac{e^{-i(\zeta-z)
k}h(\zeta)d\zeta}{\zeta-z}\,\,\, z\in \C-\Om.
\end{align}
The fact that $h$ solves \eqref{outside} implies that $\lim_{z\to
z_0\in\del\Om} W_r(z)=h(z_0)$. Thus $W_r$ is an outside continuous
extension of $h$. Moreover, $e^{-izk} W_r-1 = O(1/z)$ for $z$
large, hence $W_r\in L^\pt(\C-\Om)$.

We produced a continuous map in $\R^2$ which solves
\eqref{dbarlong} both inside and outside $\Om$ and behaves like
$e^{izk}$ for $z$ large. We need to check that it solves the
equation \eqref{dbarlong} across the boundary. Since $b$ has
compact support inside $\Om$ we have that $W_r$ is in fact
analytic in both sides of the boundary and continuous across.
Morera's theorem asserts that $W_r$ must be in fact analytic
across. Therefore $W_r$ solves \eqref{dbarlong} in the whole plane
and has the right behavior at infinity. Uniqueness in lemma
\ref{Ws} concludes the proof.
\end{proof}

Immediate consequence to the proposition \ref{relation} and to the
pointwise estimates \eqref{longKestimate} we can determine the
traces on $\del\Om$ of $\Psi_r$ and $\Psi_i$. Moreover by formula
\eqref{newScatT} we determine the scattering transform.
\begin{corollary}[Reconstruction of the scattering transform]\label{recon_t}
Under the assumptions of the proposition \ref{relation}
we have
\begin{align}\label{einvb}
e^{-\dbinv b}(z)= \lim_{k\to\infty}W_r(z,k),\qquad z\in\del\Om.
\end{align}
and for any $k\in\C$ we recover
\begin{align}
&\Psi_r(z,k)=e^{\dbinv b}(z)W_r(z,k),\qquad z\in\del\Om,\\
&\Psi_i(z,k)=e^{\dbinv b}(z)W_i(z,k),\qquad z\in\del\Om.
\end{align}Moreover,
\begin{align}
t(k)=\frac{i}{2\pi}\int_{\del\Om}e^{i\ol{zk}}\ol{\nu}(z)
\left(\ol{\Psi_r}(z,k)-i\ol{\Psi_i}(z,k)\right)d\sigma(z),
\end{align}is a function in $L^r(\R^2)\cap L^{\rt}(\R^2)\cap
L^{r'}(\R^2)$ for some $r<2$, $\rt^{-1}=r^{-1}-1/2$ and
$r^{-1}+r'^{-1}=1$.
\end{corollary}
Now we use the inverse scattering method of theorem 1.2 to
reconstruct $q$.
\begin{corollary}Let $\Phi_r\sim e^{izk}$ and $\Phi_i\sim
ie^{izk}$ in $L^\pt$ for large $k\in\C$ be the unique solutions
\begin{equation}
\frac{\del \Phi}{\del \ol{k}}(k)+t(k)\ol{\Phi}(k)=0,\,\, k\in \C.
\end{equation}Then
\begin{align}
q(z)=-\frac{i}{\pi}\int_{\R^2}e^{i\ol{zk}}\ol{t}(k)\left(
\Phi_r(z,k)-i\Phi_i(z,k)\right)d\mu(k).
\end{align}
\end{corollary}
Knowing $q$ we also know $|b|$ since from \eqref{new_q} we have
$|q|=|b|$. Next we show how to determine its argument by solving
\eqref{new_q} to recover $b$.

The following result is due to Cheng and Yamamoto
\cite{chengYamamoto98}. For the sake of completeness we sketch its
proof.
\begin{lemma}If $q\in L^\pt(\Om)$ then there exist at most one
solution $b\in L^\pt(\Om)$ of the equation
\begin{align}\label{recon_q}
\ol{b}(z)e^{\dbinv b (z) -\deinv\ol{b}(z)}=q(z),\,\, z\in \Om.
\end{align}
\end{lemma}
\begin{proof}
Assume there are two solutions $b_1,b_2\in L^\pt(\Om)$ and let
$d=\dbinv(b_2-b_1)\in W^{1,\pt}(\Om)\subset H^1(\Om)$.

From \eqref{einvb} we have $d=0$ on $\del\Om$. Hence $d\in
H^{1}_0(\Om)$. Since both solve \eqref{recon_q} we have
$\ol{b_1}(z)=\ol{b_2}(z)e^{d(z)-\ol{d}(z)}$, from where
\begin{align}
|\dba d|=|b_2|\cdot|e^{d(z)-\ol{d}(z)}-1|=|q|\cdot
|e^{d(z)-\ol{d}(z)}-1|\leq |q|\cdot|d-\ol{d}|\leq 2|q||d|.
\end{align}

By Carleman estimates for $d\in H^1(\R^2)$ of compact support (see
H\"{o}rmander \cite{hormander}, Prop. 17.2.3) we have
\begin{align}
\int_{\Om}\Delta\varphi|d|^2e^{2\varphi}dx\leq 4\int_{\Om}|\dba
d|^2e^{2\varphi}dx,
\end{align}for some strictly convex function $\varphi\in
C^2(\Om)$. Approximate a $\varphi\in H^{2,\pt}(\Om)$ solution of
$\Delta \varphi=17|q|^2$ in $L^2(\Om)$ by a smooth sequence
$\varphi_n\to\varphi$ uniformly on $\ol{\Om}$. Then
$$\int_{\Om}\Delta\varphi_n|d|^2e^{2\varphi_n}dx\leq 16\int_{\Om}|q|^2|d|^2e^{2\varphi_n}dx.$$
For $n$ sufficiently large the reverse inequality holds. Hence
$d=0$ and $b_1=b_2$.
\end{proof}
We are left to find the unique solution of \eqref{recon_q}.
\begin{lemma}[Phase unwrapping] Let $v\in 1+ L^\pt(\R^2)$ be the unique solution of
\begin{align}
\dba v =\ol{qv}
\end{align}
then $v$ vanishes on a set of measure zero. Define
$b=\ol{q}\ol{v}/v$ on the set where $v$ does not vanish, else we
can set $b=q$. Then $b$ is the unique solution of \eqref{recon_q}
in $L^\pt(\Om)$.
\end{lemma}
\begin{proof}Existence and uniqueness of $v$ follows from the Fredholm alternative as before.
It is known from Vekua \cite{vekua} that the set of zeroes of
pseudo-analytic functions has measure zero . Since $\ol{bv}=qv$ we
have that $v$ also solves $\dba v= bv$ in the whole plane.
Equivalently $\dba(e^{-\dbinv b}v)=0$. Thus $e^{-\dbinv b}v$ is
analytic and also goes to $1$ as $|z|\to\infty$. By Liouville's
theorem we have $v=e^{\dbinv b}$. From its definition we have
$$b= \ol{q}e^{-\del^{-1}\ol{b} +\dbinv b}.$$
\end{proof}
\section{Concluding Remarks}
In order to solve the inverse problem one finds first the traces
of the exponentially behaving solutions. It is easy to show that
any solution of \eqref{outside} outside a disk can be represented
as a series
$$W(z,k)=e^{izk}\sum_{n=0}^\infty \frac{a_n}{z^n},$$ with $a_n$'s unknown
coefficients. We determine them by solving the singular boundary
integral equations \eqref{inside}. This step is severely ill posed
and regularization techniques are necessary, truncation in the
above series helps, see Knudsen \cite{kim} for further ideas of
regularization. Moreover, there is only a logarithmic type
stability, see Barcelo et. al. in \cite{bbr}.

The second step consists in constructing the scattering transform
$t(k)$ via the formulae of corollary \ref{recon_t}. Next we solve
the weakly singular integral equations \eqref{freq_eq} in the
$k$-space. This part is stable. It is here that we need the
$\epsilon$-extra regularity. One needs $t\in L^r(\R^2_k)$ for some
$r<2$ in order to solve \eqref{freq_eq}. If $q$ is only in
$L^\pt_c$ then $t\in L^2(\R^2)$ (according to Sung \cite{sung2} as
corrected by Brown and Uhlmann \cite{brownUhlmann}) and this
suffices for uniqueness. This would recover the uniqueness result
of Cheng and Yamamoto. It is not clear how to find solutions of
\eqref{freq_eq} when $t\in L^2(\R^2_k)$.

Reconstruct $q$ from the formula \eqref{fourierInv}. Notice that
we have estimates of decay in $k$ for $t\in L^r (\R^2_k)$ as well
as for $e^{-izk}(\Phi_r -i\Phi_i)-2$ as given in
\eqref{extrasmooth}. These can lead to estimates of the truncation
error in the integral in \eqref{fourierInv}.

One of the questions in \cite{uhlmann} concerned the
characterization of traces of exponentially behaving solutions in
the first order system in $\Om$: $\dba v=q w$ and $\del
w=\ol{q}v$. A partial answer was given by Knudsen and the author
in \cite{knudsenTamasan} for $q$ of the special form $q=\del f$
with $f$ real valued. We can give now the answer for a general
$q$. Note that $v\pm\ol{w}$ solves the $\dba$-equation $\dba u
+q\ol{u}=0$, and that we characterized the traces on the boundary
of such solutions in terms of a Hilbert transform.
\begin{center}
{\bf Acknowledgement}
\end{center}
I would like to thank Professor A. Nachman for his generous
sharing of ideas in the inverse scattering theory, and the
organizers, Professors J. McLaughlin and H. Engl, for the
invitation in the special semester on Inverse Problems at
IPAM-UCLA.

\end{document}